\documentstyle[twoside,amssymb,12pt]{article}
\def\Q{{\Bbb Q}}

\def\Z{{\Bbb Z}}

\newtheorem{prop}{Proposition}[section]
\newtheorem{dfn}[prop]{Definition}
\newtheorem{theo}[prop]{Theorem}

\newtheorem{rem}[prop]{Remark}
\newtheorem{coro}[prop]{Corollary}

\title{\sc Stringy Hodge numbers and Virasoro algebra} 

\author{{\sc Victor V. Batyrev} \\
\small  {\em Mathematisches Institut, Universit\"at T\"ubingen}   \\
\small  {\em Auf der Morgenstelle 10,  72076  T\"ubingen, Germany}  \\
\small  {\em e-mail: batyrev@bastau.mathematik.uni-tuebingen.de} \\
 }

\begin{document}

\date{}

\maketitle

\begin{abstract}
Let $X$ be an arbitrary  smooth $n$-dimensional projective variety. 
It was discovered by Libgober and Wood that the product of the Chern
classes $c_1(X)c_{n-1}(X)$ depends only on the  
Hodge numbers of $X$. This result has been used by  Eguchi, Jinzenji 
and Xiong in their approach to the quantum cohomology of $X$ via 
a representation of the Virasoro algebra with  the central charge 
$c_n(X)$. 

In this paper we define for singular  varieties $X$ a rational number 
$c_{\rm st}^{1,n-1}(X)$ which is a stringy version of the number $c_1c_{n-1}$ 
for smooth $n$-folds. We show that the number $c_{\rm st}^{1,n-1}(X)$ 
can be expressed in the same way using the   
stringy Hodge numbers of $X$. Our results provides  an  evidence 
for the existence of an   
approach to quantum cohomology of singular varieties $X$ 
via a representation of  the Virasoro algebra whose    
central charge is the rational number 
$e_{\rm st}(X)$ which equals  the stringy Euler number of $X$. 
\end{abstract}

\thispagestyle{empty}

\newpage

\section{Introduction}

Let $X$ be an arbitrary  smooth projective variety of dimension $n$. 
The $E$-polynomial of $X$ is defined as 
\[ E(X; u,v):= \sum_{p,q} (-1)^{p+q} h^{p,q}(X) u^p v^q \]
where $h^{p,q}(X)= dim\, H^q(X, \Omega^p_X)$ are Hodge numbers of $X$. 
Using the Hirzebruch-Riemann-Roch theorem, Libgober and Wood \cite{LW} 
has proved the following equality (see also the papers of Borisov 
\cite{LB1} and Salamon \cite{S}): 

\begin{theo} 
\[ \frac{d^2}{d u^2}{E}_{\rm st}(X; u,1)|_{u =1} = \frac{3n^2 - 5n}{12} 
c_n(X) + \frac{c_1(X)c_{n-1}(X)}{6}.  \]
\label{d-rel} 
\end{theo} 

\noindent
By Poincar{\'e} duality for $X$, one immediatelly obtains \cite{LB1,S}:  

\begin{coro} 
Let $X$ be an arbitrary smooth $n$-dimensional projective variety. Then 
 $c_1(X)c_{n-1}(X)$ can be expressed  via  
the Hodge numbers of $X$ using  the following 
equality 
\[ \sum_{p,q} (-1)^{p+q} h^{p,q}(X) \left(p - \frac{n}{2} \right)^2 = 
\frac{n}{12}c_n(X)  + 
\frac{1}{6} c_1(X)c_{n-1}(X), \]
where 
\[ c_n(X) = \sum_{p,q}  (-1)^{p+q}h^{p,q}(X) \]
is the Euler number of $X$. 
\label{vir-rel}
\end{coro}

\begin{coro} 
Let $X$ be an arbitrary 
 smooth $n$-dimensional projective variety with $c_1(X) =0$. 
Then the Hodge numbers of $X$  satisfy the following 
equation
\[ \sum_{p,q} (-1)^{p+q} h^{p,q}(X) \left(p - \frac{n}{2} \right)^2 = 
\frac{n}{12} \sum_{p,q}  (-1)^{p+q}h^{p,q}(X),  \]
In particular, for hyper-K\"ahler manifolds $X$ this equation reduces to 
\[ 2 \sum_{j =1}^{2n} (-1)^j (3j^2 - n) b_{2n-j}(X) = nb_{2n}(X),
\]
where 
\[ b_i(X) =  \sum_{p+q=i}  h^{p,q}(X) \]
is $i$-th Betti number of $X$.   
\label{cy-rel}
\end{coro}

\begin{rem} 
{\rm We that if  $X$ is a $K3$-surface, then  the  relation \ref{cy-rel} 
is equivalent to the equality $c_2(X) =24$. 
For smooth 
Calabi-Yau $4$-folds $X$ the relation \ref{cy-rel} has been observed 
by Sethi,  Vafa, and  Witten \cite{W2} (it is  equivalent to the equality 
\[  c_4(X) = 6( 8 - h^{1,1}(X) + h^{2,1}(X) - h^{3,1}(X)), \]
if  $h^{1,0}(X) = h^{2,0}(X) = 
h^{3,0}(X) =0$). }  
\end{rem}

There are a lot of examples of Calabi-Yau varieties $X$ having at worst 
 Gorenstein canonical singularities which are hypersurfaces and 
complete intersections in Gorenstein 
toric Fano varieties \cite{BA,BB0}. 
It has been shown in \cite{BD} that for all these examples 
of singular Calabi-Yau varieties $X$ one can define so called {\em stringy 
Hodge numbers} $h^{p,q}_{\rm st}(X)$. 
Moreover,  the stringy Hodge numbers of Calabi-Yau 
complete intersections in Gorenstein toric varieties agree with 
the topological mirror duality test \cite{BB}. 
It was a natural question posed in 
\cite{B1}, whether one has the same identity for stringy Hodge numbers 
of singular Calabi-Yau varieties as for usual Hodge numbers of smooth 
Calabi-Yau manifolds, i.e. 

\begin{equation} 
 \sum_{p,q} (-1)^{p+q} h^{p,q}_{\rm st}(X)
 \left(p - \frac{n}{2} \right)^2 = 
\frac{n}{12} \sum_{p,q}  (-1)^{p+q}h^{p,q}_{\rm st}(X) = \frac{n}{12} 
e_{\rm st}(X).
\label{st-rel}
\end{equation} 

\noindent
The purpose of this paper is to show that the formula 
(\ref{st-rel}) holds true. 
Moreover, one can define a rational number $e_{\rm st}^{1, n-1}(X)$ 
which is a stringy version $c_1(X)c_{n-1}(X)$ such that the stringy 
analog of \ref{vir-rel}

\begin{equation}
\sum_{p,q} (-1)^{p+q} h^{p,q}_{\rm st}(X) \left(p - \frac{n}{2} 
\right)^2 = 
\frac{n}{12}e_{\rm st}(X)  + 
\frac{1}{6} c_{\rm st}^{1,n-1}(X) 
\label{st-rel2}
\end{equation}
holds true provided the  stringy Hodge numbers of $X$ exist.

\section{Stringy Hodge numbers} 

Recall our general approach to the notion of stringy Hodge numbers 
$h^{p,q}_{\rm st}(X)$ for projective algebraic varieties $X$ with 
canonical singularites (see \cite{B1}).  Our main definition in 
\cite{B1} can be reformulated as follows: 

\begin{dfn} 
{\rm Let $X$ be an arbitrary $n$-dimensional projective 
variety with at worst log-terminal singularites, 
$\rho\, : \, Y \to X$ a resolution of singularities whose exceptional 
locus $D$ is a divisors with normally crossing components 
$D_1, \ldots, D_r$. We set $I: = \{1, \ldots, r \}$ and 
$D_J: = \bigcap_{j \in J} D_j$ for all $J \subset I$. 
Define the {\bf stringy $E$-function} of $X$ to be 
\[ E_{\rm st} (X; u,v) := 
\sum_{J \subset I} E(D_J; u,v) \prod_{j \in J}
\left(\frac{uv-1}{(uv)^{a_j+1} -1} -1 \right), \]
where the rational numbers $a_1, \ldots, a_r$ are determined by the equality
\[ K_Y = \rho^*K_X + \sum_{i =1}^r a_i D_i. \]
Then the {\bf stringy Euler number} of $X$ is defined as 
\[ e_{\rm st}(X) := 
\lim_{u,v \to 1} E_{\rm st} (X; u,v) = \sum_{J \subset I} 
c_{n - |J|}(D_J) \prod_{j \in J}
\left(\frac{-a_j}{a_j+1} \right), \]
where $c_{n-|J|}(D_J)$ is the Euler number of $D_J$ (we set 
$ c_{n-|J|}(D_J)=0$ if $D_J$ is empty).   
}
\label{def-e}
\end{dfn} 

\begin{dfn} 
{\rm Let $X$ be an arbitrary $n$-dimensional projective  
variety with at worst Gorenstein canonical singularites. 
We say that {\bf stringy Hodge numbers of $X$ exist}, if 
$E_{\rm st}(X; u,v)$ is a polynomial, i.e., 
\[ E_{\rm st}(X; u,v) = \sum_{p,q} a_{p,q}(X)u^pv^q .\] 
Under the assumption that $E_{\rm st}(X; u,v)$ is a polynomial, we 
define the {\bf stringy Hodge numbers} $h^{p,q}_{\rm st}(X)$ to be 
$(-1)^{p+q}a_{p,q}$. 
} 
\end{dfn}

\begin{rem} 
{\rm In the above definitions, the condition that $X$ has at worst 
log-terminal singularities means that $a_i > -1$ for all $i \in I$; 
the condition that  $X$ has at worst 
 Gorenstein canonical  singularities is equivalent for 
$a_i$ to be nonnegative integers  for all $i \in I$ (see \cite{KMM}).}
\end{rem}

The following statement has been proved in \cite{B1}: 

\begin{theo} 
Let $X$ be an arbitrary $n$-dimensional projective  
variety with at worst Gorenstein canonical singularites. 
Assume that  stringy Hodge numbers of $X$ exist. 
Then they have the following properties: 

{\rm (i)} $h^{0,0}_{\rm st}(X)= h^{n,n}_{\rm st}(X) =1$; 

{\rm (ii)} $h^{p,q}_{\rm st}(X)=h^{n-p,n-q}_{\rm st}(X)$ and 
 $h^{p,q}_{\rm st}(X)=h^{q,p}_{\rm st}(X)$ $\forall p,q$; 

{\rm (iii)} $h^{p,q}_{\rm st}(X) = 0$ $\forall p,q >n$.
\label{st-pr} 
\end{theo}

\section{The number $c_{\rm st}^{1,n-1}(X)$}

\begin{dfn} 
{\rm Let $X$ be an arbitrary $n$-dimensional projective variety $X$ having 
 at worst log-terminal singularities and $\rho\;:\; Y \to X$ is a 
desingularization with normally crossing irreducible components 
$D_1,\ldots, D_r$ of the exceptional locus.   
We define the number 
\[  c_{\rm st}^{1, n-1}(X): =  \sum_{J \subset I} \rho^*c_1(X) 
c_{n-|J|-1}(D_J)  
\prod_{j \in J} \left( \frac{-a_j}{a_j +1} \right), \]
where  $\rho^*c_1(X) 
c_{n-|J|-1}(D_J)$ is considered as the intersection number of the $1$-cycle 
$c_{n-|J|-1}(D_J) \in A_1(D_J)$ with the $\rho$-pullback of 
the class of the anticanonical $\Q$-divisor of $X$.}
\label{def-c1}    
\end{dfn}  

\begin{rem} 
{\rm It is not clear a priori that  
the number $c_{\rm st}^{1, n-1}(X)$ in the above the definition  does not
depend on the choice of a desingularization $\rho$. Later we shall see that 
it is the case.} 
\end{rem}  

\begin{prop} 
For any smooth $n$-dimensional projective variety $V$, one has 
\[ \frac{d}{d u} {E}(V;u,1)|_{u =1} = \frac{n}{2} 
c_n(V). \]
\label{f-1}
\end{prop} 

\noindent
{\em Proof.} By definition of $E$-polynomials, we have 
\[  \frac{d}{d u} {E}(V;u,1)|_{u =1} = 
\sum_{p,q} p(-1)^{p+q}h^{p,q}(V). \]
The  Poincar{\'e} duality $h^{p,q}(V) = h^{n-p,n-q}(V)$ $\forall p,q$ 
implies that 
\[ \sum_{p,q} \left(p - \frac{n}{2} \right)  (-1)^{p+q}h^{p,q}(V) = 0. \] 
Hence, 
\[ \sum_{p,q} p(-1)^{p+q}h^{p,q}(V) = \frac{n}{2} \sum_{p,q}
 (-1)^{p+q}h^{p,q}(V) =\frac{n}{2} c_n(V). \] 
\hfill $\Box$ 

\begin{prop} 
For any $n$-dimensional projective variety $X$ having 
 at worst log-terminal singularities, one has 
\[ \frac{d}{d u} {E}_{\rm st}(X;u,1)|_{u =1} = \frac{n}{2} 
e_{\rm st}(X). \]
\label{f-2}
\end{prop} 

\noindent
{\em Proof.} By definition \ref{def-e}, we have 
   \[ E_{\rm st} (X; u,1) = \sum_{J \subset I} E(D_J; u,1) 
\prod_{j \in J} \left( \frac{u -1}{u^{a_j +1} -1} -1 
\right). \] 
Applying \ref{f-1} to every  smooth submanifold $D_J \subset Y$, we obtain
\[ \frac{d}{d u} {E}_{\rm st}(X;u, 1)|_{u =1} = 
\sum_{J \subset I} \frac{(n-|J|)}{2} c_{n-|J|}(D_J) 
\prod_{j \in J} \left( \frac{-a_j}{a_j +1} \right)  + \] 
\[ +  \sum_{J \subset I}  \frac{|J|}{2} c_{n-|J|}(D_J) \prod_{j \in J} 
\left( \frac{-a_j}{(a_j +1)} \right) = 
 \frac{n}{2}  \sum_{J \subset I} c_{n-|J|}(D_J) 
\prod_{j \in J} \left( \frac{-a_j}{a_j +1} \right) = 
\frac{n}{2}e_{\rm st}(X). \]

\hfill $\Box$ 

\begin{prop} 
Let $V$ be a smooth projective algebraic variety of 
dimension  $n$ and $W \subset V$ a smooth irreducible 
divisor on $V$ or empty divisor $($the latter means that 
${\cal O}_V(W) \cong {\cal O}_V)$. Then 
\[ c_1({\cal O}_V(W)) c_{n-1}(V) = c_{n-1}(W) +  
c_1({\cal O}_W(W)) c_{n-2}(W), \]
where  $c_{n-1}(W)$ is considered to be zero if $W = \emptyset$. 
\label{div}
\end{prop}

\noindent
{\em Proof.} Consider the short exact sequnce 
\[ 0 \to T_W \to T_V|_W \to  {\cal O}_W(W) \to 0, \]
where $T_W$ and $T_V$ are tangent shaves  on $W$ and $V$.  
It gives the following the relation betwen Chern polynomials 
\[ (1 + c_1({\cal O}_W(W) t) ( 1 + c_1(D)t + c_2(D)t^2 + \cdots + 
c_{n-1}(D)t^{n-1}) = \]
\[ = 1 + c_1( T_V|_W)t + c_2( T_V|_W)t^2 + 
c_{n-1}( T_V|_W) t^{n-1} ). \]
Comparing the coefficients by $t^{n-1}$ and using 
$c_{n-1}( T_V|_W) =  c_1({\cal O}_V(W)) c_{n-1}(V)$, we come 
to the required equality. 
\hfill $\Box$ 

\begin{coro}
Let $Y$ be a smooth projective variety, $D_1, \ldots, D_r$ 
smooth irreducible divisors with normal crossings, 
$I: = \{1, \ldots, r \}$. Then for all  $J \subset I$ and for 
all $j \in J$ 
one has 
\[ c_1({\cal O}_{D_{J 
\setminus \{j \} }}(D_j)) c_{n-|J|}( D_{J \setminus \{j \} }) - 
 c_{n-|J|}( D_{J})  = c_1({\cal O}_{D_J}(D_j))c_{n-|J|-1}(D_J), \]
where  
$D_J$ is the complete intersection $\bigcap_{j \in J} D_j$. 
\label{rec} 
\end{coro} 

\noindent
{\em Proof.} One sets in \ref{div} $V:=  D_{J \setminus \{j \}}$ and 
$W:=  D_{J}$. 
\hfill $\Box$

\begin{prop} 
Let $\rho\; : \; Y \to X$ be a desingularization as in \ref{def-c1}. Then
\[  \sum_{J \subset I} c_1(D_J) c_{n-|J|-1}(D_J)  
\prod_{j \in J} \left( \frac{-a_j}{a_j +1} \right) = c_{\rm st}^{1, n-1}(X) +\]
\[ + 
 \sum_{J \subset I} \left(  \sum_{j \in J}  (a_{j} +1)c_{n-|J|}(D_J) 
 \right)\prod_{j \in J} \left( \frac{-a_j}{a_j +1} \right) . \]
\label{relat-div}
\end{prop} 

\noindent
{\em Proof.} Using  the formula 
\[ c_1(Y) = \rho^*c_1(X) + \sum_{i \in I} -a_i c_1({\cal O}_Y(D_i)) \]
and the adjunction formula for every  complete intersection 
$D_J$ $(J \subset I)$, we obtain 
\[ c_1(D_J) = \rho^*c_1(X)|_{D_J} + 
\sum_{j \in J} (-a_j-1)c_1({\cal O}_{D_J}(D_j)) + 
\sum_{j \in I \setminus J} (-a_j)c_1({\cal O}_{D_J}(D_j)). \]
Therefore 
\begin{equation}
  \sum_{J \subset I} c_1(D_J) c_{n-|J|-1}(D_J)  
\prod_{j \in J} \left( \frac{-a_j}{a_j +1} \right) = c_{\rm st}^{1, n-1}(X) +
\label{eq3}
\end{equation}
\[ + \left( \sum_{j \in J} (-a_j-1)c_1({\cal O}_{D_J}(D_j)) 
c_{n-|J|-1}(D_J) \right) \prod_{j \in J} \left( \frac{-a_j}{a_j +1} 
\right)  +\]
\[ +  \left( \sum_{j \in I \setminus J} (-a_j)c_1({\cal O}_{D_J}(D_j))
c_{n-|J|-1}(D_J) \right)
\prod_{j \in J} \left( \frac{-a_j}{a_j +1} \right) . \]
Using \ref{rec}, we obtain 
\begin{equation}
  \sum_{j \in J} (-a_j-1)c_1({\cal O}_{D_J}(D_j)) 
c_{n-|J|-1}(D_J)  = 
\label{eq4}
\end{equation}
\[ =   \sum_{j \in J} (-a_j-1) \left(c_1({\cal O}_{D_{J 
\setminus \{j \} }}(D_j)) c_{n-|J|}( D_{J \setminus \{j \} }) - 
 c_{n-|J|}( D_{J}) \right). \]
By substitution (\ref{eq4}) to (\ref{eq3}), we come to the 
required equality. 
\hfill $\Box$

\begin{theo} 
Let  $X$ be  an arbitrary 
$n$-dimensional projective variety variety with 
at worst log-terminal singularities.  Then 
\[ \frac{d^2}{d u^2}{E}_{\rm st}(X; u,1)|_{u =1} = \frac{3n^2 - 5n}{12} 
e_{\rm st}(X) + \frac{1}{6}c_{\rm st}^{1,n}(X). \]
\label{main}
\end{theo} 

\noindent
{\em Proof.} Using the equalities 
\[  \frac{d}{d u}\left( \frac{u-1}{u^{a+1} -1} -1 \right)_{u=1} = 
\frac{-a}{2(a+1)}, \;\; 
 \frac{d^2}{d u^2} \left(  \frac{u-1}{u^{a+1} -1} -1 
\right)_{u=1} = \frac{a(a+2)}{6(a+1)}  \]
together with the identities in \ref{d-rel} and 
\ref{f-1} for every submanifold 
$D_J \subset Y$, we obtain

\[ \frac{d^2}{d u^2} {E}_{\rm st}(X; u,1)|_{u =1} = 
 \sum_{J \subset I}  \frac{c_1(D_J)c_{n - |J|-1}(D_J)}{6} 
\prod_{j \in J} \left( \frac{-a_j}{a_j +1} \right) + \]
\[ + c_{n-|J|}(D_J) \frac{3(n-|J|)^2 - 5(n-|J|)}{12} 
\prod_{j \in J} \left( \frac{-a_j}{a_j +1} \right) + \]
\[ +   \sum_{J \subset I}  
 \frac{(n-|J|)|J|c_{n-|J|}(D_J)}{2} \prod_{j \in J} 
\left( \frac{-a_j}{a_j +1} \right) + \]
\[ +  \sum_{J \subset I} \frac{c_{n-|J|}(D_J)(|J|-1)|J|}{4} 
\prod_{j \in J} 
\left( \frac{-a_j}{a_j +1} \right) + \]
\[ + 
\sum_{J \subset I} \frac{c_{n-|J|}(D_J)(-\sum_{j \in J} (a_j +2) )}{6} 
\prod_{j \in J} 
\left(  \frac{-a_j}{a_j+1}  \right). \]
By \ref{relat-div}, the first term of the above  equals  
\[  
\frac{1}{6} c_{\rm st}^{1, n-1}(X) +   \frac{1}{6} 
 \sum_{J \subset I} \left(  \sum_{j \in J}  (a_{j} +1)c_{n-|J|}(D_J) 
\right)\prod_{j \in J} \left( \frac{-a_j}{a_j +1} \right). \]
Now the required statement  follows from the equality
\[ \frac{  \sum_{j \in J}  (a_{j} +1)}{6} 
  +   \frac{3(n-|J|)^2 - 5(n-|J|)}{12} +  
\frac{(n-|J|)|J|}{2} + \]
\[ +  \frac{(|J|-1)|J|}{4} +
\frac{-\sum_{j \in J} (a_j +2) }{6} = \frac{3n^2 - 5n}{12}. 
   \]
\hfill $\Box$ 

\begin{coro} 
The number $c_{\rm st}^{1,n}(X)$ does not depend on the choice of 
the desingularization $\rho\;: \; Y \to X$. 
\end{coro} 

\noindent
{\em Proof.} By \ref{f-2} and \ref{main},  $c_{\rm st}^{1,n}(X)$ 
can be computed in terms of derivatives of the stringy $E$-function 
of $X$. But the stringy $E$-function does not depend on the choice 
of a desingularization \cite{B1}. 
\hfill $\Box$

\begin{coro} 
Let $X$ be a projective variety with at worst Gorenstein canonical 
singularities. Assume that the stringy Hodge numbers of $X$ exist. 
Then 
\[ \sum_{p,q} (-1)^{p+q} h^{p,q}_{\rm st}(X) \left(p - \frac{n}{2} 
\right)^2 = 
\frac{n}{12}e_{\rm st}(X)  + 
\frac{1}{6} c_{\rm st}^{1,n-1}(X). \]
\end{coro} 

\noindent
{\em Proof.} The equality follows immediately from \ref{main} 
using the properties of the stringy Hodge numbers \ref{st-pr}.
\hfill $\Box$ 

\begin{coro}
If the canonical class of $X$ is numerically trivial, then 
 $c_{\rm st}^{1,n}(X) =0$. In particular, for Calabi-Yau varieties with 
at worst Gorenstein canonical singularities we have 
\[ \frac{d^2}{d u^2}{E}_{\rm st}(X; u,1)|_{u =1} = \frac{3n^2 - 5n}{12} 
e_{\rm st}(X), \] 
and therefore stringy Hodge numbers of $X$ satisfy the identity   
$($\ref{st-rel}$)$ provided these stringy numbers  exist. 
\end{coro}

\section{Virasoro Algebra} 

Recall that the Virasoro algebra with the central charge $c$ 
consists of operators $L_n$ $(m \in \Z)$ satisfying the 
relations 
\[ [L_n, L_m] = (n-m)L_{n+m} + c\frac{n^3-n}{12}\delta_{n+m,0} 
\;\;n, m \in \Z. \] 

For arbitrary compact  K\"ahler manifold $X$, 
Eguchi et. al have proposed in \cite{EHX,EMX} a new approach 
to its quantum cohomology and to its Gromov-Witten invariants for 
all genera $g$ using so called the Virasoro condition: 
\[ L_n Z = 0, \forall n \geq -1, \]
where 
\[ Z = {\rm exp}\, F =  {\rm exp}\left(
 \sum_{g \geq 0} \lambda^{2g-2}F_g \right) \]
 is the partition function of the topological $\sigma$-model
with the target space $X$ and $F_g$ the free energy function corresponding 
to the genus $g$. In this approach, the central charge 
$c$ acts as the multiplication by $c_n(X)$. Moreover, 
all  Virasoro operators $L_n$ can be 
explicitly written in terms of elements of a  
basis of the cohomology of $X$, their  gravitational descendants and 
the action of  $c_1(X)$ on the cohomology by the multiplication. 
In particular the commutator relation  
\[ [L_1, L_{-1}] = 2 L_0 \]
implies the precisely the identity of Libgober and Wood in the form 
\[ \sum_{p,q}(-1)^{p+q} h^{p,q}(X) 
 \left( \frac{n+1}{2} - p \right) \left(p - \frac{n-1}{2} \right) = 
\frac{1}{6} \left( \frac{3-n}{2} c_n(X) - c_1(X)c_{n-1}(X) \right).  \]

Now let $X$ be a projective algebraic variety with at worst log-terminal 
singularities. We conjecture  that there exists an analogous approach 
to the quantum cohomology as well as  to the  
Gromov-Witten invariants of  $X$ for all genera using the  Virasoro algebra  
in such a way  that for any resolution of singularities 
$\rho\; : \; Y \to X$ the corresponding  
Virasoro operators can be explicitely computed 
via the   numbers $a_i$ appearing in the formula 
\[ K_X = \rho^*K_X + \sum_{i =1}^r a_i D_i \]
and bases in cohomology of all complete intersections $D_J$ together
with the multiplicative actions of $c_1(D_J)$ in them.  
We consider our main result \ref{main} as an evidence in favor of this 
conjecture.

\end{document}